# TESTING FOR MONOTONE INCREASING HAZARD RATE

By Peter Hall and Ingrid Van Keilegom[1]

*Australian National University, and Australian National University and Université Catholique de Louvain*

A test of the null hypothesis that a hazard rate is monotone nondecreasing, versus the alternative that it is not, is proposed. Both the test statistic and the means of calibrating it are new. Unlike previous approaches, neither is based on the assumption that the null distribution is exponential. Instead, empirical information is used to effectively identify and eliminate from further consideration parts of the line where the hazard rate is clearly increasing; and to confine subsequent attention only to those parts that remain. This produces a test with greater apparent power, without the excessive conservatism of exponential-based tests. Our approach to calibration borrows from ideas used in certain tests for unimodality of a density, in that a bandwidth is increased until a distribution with the desired properties is obtained. However, the test statistic does not involve any smoothing, and is, in fact, based directly on an assessment of convexity of the distribution function, using the conventional empirical distribution. The test is shown to have optimal power properties in difficult cases, where it is called upon to detect a small departure, in the form of a bump, from monotonicity. More general theoretical properties of the test and its numerical performance are explored.

**1. Introduction.** Estimation of a hazard rate under the hypothesis that it is nondecreasing, and testing the validity of this assumption, are motivated by problems where failure rate of a machine part or a biological system can be expected to increase with lifetime. If for some reason a machine part becomes more reliable with time over at least part of its life cycle, then it can be particularly important to know that fact. The knowledge may lead to changes in the way the part is manufactured or finished, so as to remove

Received September 2002; revised May 2004.
[1]Supported in part by the contract "Projet d'Actions de Recherche Concertées" nr 98/03-217 and by the IAP research network nr P5/24 of the Belgian state.

the requirement for a running-in period where failure is relatively likely to occur. In this paper we suggest a new test statistic of the null hypothesis of monotone nondecreasing failure rate and a new approach to calibrating the distribution of the statistic so as to determine a critical point for the test.

Our methods confer two advantages relative to existing approaches. First, our test statistic is focused on relatively "local" departures from the null hypothesis of nondecreasing hazard rate, and pays relatively little attention to those parts of the sample space where the hazard rate is indeed monotone nondecreasing. Nevertheless, the method is easily localized still further, since it focuses on variation of the hazard rate over an interval which can be increased or decreased at the investigator's discretion, or, indeed, replaced by the union of two or more intervals.

Second, our new method of calibration makes the test statistic much more sensitive to relatively small departures from the null hypothesis. For a given nominal probability of rejection, our calibration approach produces a test with greater apparent power than do standard methods based on calibration by comparison with the exponential distribution. The reason is that the exponential case is particularly awkward to detect; the corresponding hazard rate is perfectly flat, and, therefore, to avoid incorrectly rejecting the null hypothesis in this case, the test statistic has to satisfy itself that there are no significant bumps on a perfectly flat line. In consequence, the test tends to overlook small bumps, for fear of committing a Type I error, and so has relatively low power against hazard rates that are nondecreasing except for small bumps.

The test we propose has substantially greater apparent power in so-called "difficult cases" (cf. [7]) than does, for example, Proschan and Pyke's [19] test, calibrated using the exponential distribution. Indeed, we shall prove that our method has optimal power in this setting. That is, it is able to detect a very small perturbation of the empirical distribution, placed at a point where it produces a small nonmonotone bump in the hazard rate, and so small that even a likelihood ratio test (requiring knowledge of the shape of the bump) is barely able to detect the bump.

Our calibration method is related to the "increasing bandwidth" approach first suggested by Silverman [20] in the case of density estimation, and used in a range of other settings since; see [6] for an application in the setting of monotone nonparametric regression. However, quite unlike those applications, we increase the bandwidth only for the purpose of calibrating the test. Our test statistic does not involve any smoothing at all and is based directly on the standard empirical distribution function.

Contributions to the problem of testing for a *constant* hazard rate against a monotone alternative include those of Bickel and Doksum [5], based, like the method of Proschan and Pyke [19], on normalized spacings; Bickel [4], on the existence of asymptotically most powerful tests; Barlow and Doksum [3], on the more general problem of testing for convex orderings; and Ahmad [1], Gail and Gastwirth [11, 12] and Klefsjö [16], who proposed tests of the hypothesis of an exponential distribution. However, these approaches share the drawbacks noted above for exponential-based methods. A related difficulty arises in the context of testing for unimodality of a probability density by calibrating against the most difficult case of a uniform density; see, for example, [15]. Hall, Huang, Gifford and Gijbels [14] have suggested methods for estimating a hazard rate under the assumption of monotonicity and surveyed earlier work on the topic.

Although our focus is on testing the null hypothesis of a monotone nondecreasing hazard rate, the case where the null asserts a monotone nonincreasing rate is related. In the former case, the smoothed empirical hazard rate estimator is guaranteed to be monotone nondecreasing for all sufficiently large bandwidths, and this property is not available in the latter setting. The property makes it particularly easy to propose a bandwidth selection rule that ensures resampling from a distribution that satisfies the null; we may start with any conventional bandwidth selector, for example, based on a plug-in rule, and steadily increase the bandwidth until the smoothed empirical distribution has a monotone nondecreasing hazard rate in the region where the test is to be conducted.

There is also a simple rule in the case where $H_0$ stipulates that the hazard rate is nonincreasing: starting with any conventional bandwidth selector, increase the bandwidth until a monotone nonincreasing hazard rate is obtained; or, if that does not occur no matter how large the bandwidth, reject the null hypothesis at this point without passing to a further step. This rule is justified by the fact that, if the hazard rate is nonincreasing, then the probability that there exists a finite bandwidth (of larger order than the conventional $n^{-1/5}$), such that the smoothed empirical hazard rate is nonincreasing, generally converges to 1 as sample size increases. Nevertheless, in the remainder of the paper we shall address only the more practically important case where $H_0$ asserts a nondecreasing hazard rate.

## 2. Methodology.

2.1. *Test statistic.* Suppose the random sample $\mathcal{X} = \{X_1, \ldots, X_n\}$ is drawn from a distribution with distribution function $F$. The standard empirical distribution function is $\widehat{F}(x) = n^{-1} \sum_i I(X_i \leq x)$, where $I(\mathcal{E})$ denotes the indicator function of an event $\mathcal{E}$. The null hypothesis that $F$ has monotone hazard rate on an interval $\mathcal{I}$ is equivalent to $H = -\log(1 - F)$ being convex on $\mathcal{I}$, and, hence, provided $F$ is twice differentiable with a nonvanishing first derivative on $\mathcal{I}$, to $H''$ being nonnegative on $\mathcal{I}$. The function $H$ is the cumulative hazard rate. Its derivative is the hazard rate.

The empirical form of $H$, $\widehat{H} = -\log(1 - \widehat{F})$, is not differentiable, however. Therefore, it makes little sense to test the null hypothesis by checking for nonnegativity of the second derivative of $\widehat{H}$. We could investigate methods based directly on smoothed forms of $\widehat{H}$, but this would not necessarily lead to tests that have good power properties; see Section 2.3. Instead we note that convexity of $H$ on $\mathfrak{I}$ is equivalent to nonnegativity of $H(x+y) + H(x-y) - 2H(x)$ for all $x$ and $y$ such that both $x+y$ and $x-y$ are elements of $\mathfrak{I}$. It is not essential to take $\mathfrak{I}$ to be an interval; it can be replaced by a disjoint union of intervals, for example. In the latter case it is, however, necessary to integrate in $T$ [defined in (2.2)] over the pairs $(x, y)$ belonging to $\mathfrak{I}$ so that $x+y$ and $x-y$ lie in the same interval as $x$.

Therefore, a test of the hypothesis of increasing hazard rate or, equivalently, of

$$(2.1) \qquad H_0 : H \text{ is convex on } \mathfrak{I},$$

is to reject $H_0$ in favor of its complement if the value of

$$(2.2) \quad T = \iint_{x,y\,:\,x+y,x-y \in \mathfrak{I}} \max\{0, 2\widehat{H}(x) - \widehat{H}(x+y) - \widehat{H}(x-y)\}^r w(x,y)\,dx\,dy$$

is "too large." The exponent $r$ is an arbitrary positive number and $w$ is a nonnegative weight function. By taking the maximum in the argument of the integral at (2.2), we have largely restricted attention to places where the sampled distribution has a decreasing hazard rate. (Here and below we use the words "increasing" and "decreasing" to mean "nondecreasing" and "nonincreasing," resp.) Further restriction will be made through our method for calibration, which uses the data to determine where the hazard rate is more likely to be increasing or decreasing.

2.2. *Calibration.* Our approach to calibration will be based on bootstrap sampling from the distribution determined by a kernel density estimator,

$$\tilde{f}(x|h) = (nh)^{-1} \sum_{i=1}^{n} K\left(\frac{x - X_i}{h}\right),$$

where $K$ is a kernel and $h$ a bandwidth. We shall choose $K$ to be a smooth, symmetric density function, its graph being of conventional bell shape. Let $\widetilde{F}$ denote the distribution function corresponding to the density $\tilde{f}$, and let $\widetilde{H} = -\log(1 - \widetilde{F})$ be the associated cumulative hazard function. Then

$$(2.3) \quad \widetilde{H}''(x) = -(d/dx)^2 \log\{1 - \widetilde{F}(x)\} = \frac{\{1 - \widetilde{F}(x)\}\tilde{f}'(x) + \tilde{f}(x)^2}{\{1 - \widetilde{F}(x)\}^2}.$$

We shall write $\widetilde{H}''(x)$ as $\widetilde{H}''(x|h)$ when it is necessary to indicate dependence on bandwidth, and, as at (2.3), we shall drop the notation $h$ from quantities such as $\tilde{f}(\cdot|h)$ when it is not necessary for our argument. An empirical approach to bandwidth choice will be employed, as follows.

Let $\hat{h}$ denote a conventional empirical bandwidth, the asymptotic size of which is $n^{-1/5}$. We shall call $\hat{h}$ the "starting bandwidth." Examples include the bandwidths selected by the bootstrap, cross-validation or plug-in methods. Steadily increase the bandwidth, starting from $\hat{h}$, and stopping on the first occasion on which $\widetilde{H}''$ does not change sign on $\mathcal{I}$. Define

$$(2.4) \quad \hat{h}_{\mathrm{crit}} = \inf\{h \geq \hat{h}: \text{ the equation } \widetilde{H}''(\cdot|h) = 0 \text{ has no solution on } \mathcal{I}\}.$$

We claim that if $\mathcal{I}$ is a compact interval, then for all sufficiently large $h$, $\widetilde{H}''(\cdot|h) > 0$ on $\mathcal{I}$, and so the set at (2.4) is not empty. Therefore, $\hat{h}_{\mathrm{crit}}$ is well defined.

To verify the claim, assume $K$ has two continuous derivatives in a neighborhood of the origin, $K(0) > 0$ and $K'(0) = 0$, and observe that as $h \to \infty$, $\tilde{f}(x|h) = h^{-1} K(0) + o_p(h^{-1})$ and $\tilde{f}'(x|h) = h^{-3} K''(0) n^{-1} \sum_i (x - X_i) + o_p(h^{-3})$, where both relations hold uniformly in $x \in \mathcal{I}$. It follows that, for all sufficiently large $h$, $\tilde{f}(x)^2 > |\tilde{f}'(x)|$ for all $x \in \mathcal{I}$. The claim that $\widetilde{H}''(\cdot|h) > 0$ on $\mathcal{I}$ now follows from (2.3).

Having computed $\hat{h}_{\mathrm{crit}}$, we repeatedly create samples of size $n$ by sampling randomly, with replacement, from the distribution with density $\tilde{f}(\cdot|\hat{h}_{\mathrm{crit}})$, and thereby repeatedly compute bootstrap values, $T^*$ say, of the statistic $T$. Arguing thus, and given a nominal probability of rejection, $\alpha$ say, for the test, we may compute a critical point $\hat{c}(a)$ defined by

$$P\{T^* > \hat{c}(\alpha)|\mathcal{X}\} = \alpha.$$

The test takes the form: reject the null hypothesis if $T > \hat{c}(\alpha)$.

2.3. *The road not taken*: *tests based on* $\widetilde{H}''$. In the test described in Sections 2.2 and 2.3 we have used smoothing methods only for calibration, not to construct the test statistic itself. An alternative approach would be to base a test directly on the property that, when $H$ is twice continuously differentiable, the null hypothesis is satisfied if and only if $H'' \geq 0$ on $\mathcal{I}$. In particular, we could construct a smoothed version, $\widetilde{H}$ say, of $\widehat{H}$ with the property that $\widetilde{H}''$ is a consistent estimator of $H''$, and reject $H_0$ if (e.g.) $S = \int_{\mathcal{I}} \{\max(0, -\widetilde{H}'')\}^2$ is "too large."

This approach has drawbacks, however. First, it requires a bandwidth to be chosen when constructing the test statistic $S$; a second bandwidth would be needed when calibrating the test, if calibration were to involve sampling from a smoothed distribution. Second, the power of the test will depend intimately on choice of the first bandwidth. Indeed, the minimum distance from the null hypothesis at which local alternative distributions can be detected by the test will generally be proportional to $n^{-1/2} h^{-c}$, where $h$ is the bandwidth employed when constructing $S$, and $c > 0$ depends on the smoothing method used. Examples of this behavior in more conventional testing problems may be found in the work of Anderson, Hall and Titterington [2], Lavergne and Vuong [18] and Delecroix, Hall and Roget [10].

### 3. Theoretical properties.

3.1. *Summary of properties.* Section 3.2 shows that, if $H$ is in the class $H_{01}$ of hazard rates for which $H''$ is bounded above zero on $\mathcal{I}$ [see (3.1)], then the statistic $T$ is of size $n^{-1}$ and asymptotically normally distributed. The bootstrap accurately captures this distribution. As the convexity of $H$ becomes more marginal, the stochastic fluctuations of $T$ increase. Thus, if $H$ is in the class $H_{02}$ [see (3.6)] of hazard rates for which $H''$ vanishes at just a finite number of discrete points in $\mathcal{I}$, then the size of $T$ increases to $O(n^{-6/7})$, and its distribution becomes nonnormal (see Section 3.3). The size of $T$ increases still further, to $O_p(n^{-1/2})$, if $H''$ vanishes on an interval, and, in particular, if $F$ is an exponential distribution. (See Section 3.7, and see the third paragraph of Section 1 for an intuitive account of difficulties experienced calibrating against the exponential distribution.) Properties of our calibration method, when $H$ is in $H_{02}$, are treated in Section 3.4, where it is shown that the asymptotic probability of rejection is bounded away from zero. (Section 4 reports numerical properties in this case.) By way of contrast, if calibration is made against the exponential distribution then, when $H$ is in $H_{01}$ or $H_{02}$, the rejection probability converges to zero (Section 3.7), implying that this approach gives ultra conservatism. Optimality of our approach for identifying small, nonmonotone "wiggles" in the hazard rate is proved in Section 3.5. The ability of our calibration method to identify a fixed departure from the null hypothesis is shown in Section 3.6.

3.2. *Strict monotonicity of hazard rate.* Throughout Section 3 we shall define the statistic $T$ by taking $r = 1$ and $w \equiv 1$ in the definition at (2.2). Let $H_{01}$ be the following subset of the class of cumulative hazard functions for which $H_0$, defined at (2.1), holds:

(3.1) $H_{01} = \{H : H'' \text{ has two continuous derivatives on } \mathcal{I} \text{ and } H'' > 0 \text{ on } \mathcal{I}\}$.

(We would mention that neither $H_{01}$ nor $H_{02}$, the latter introduced in Section 3.3, is closed.) Put $g = f^{1/2}/(1 - F)$,

$$\mu = -\int_{\mathcal{I}} dx \int_{-\infty}^{\infty} E[\min\{0, y^2 H''(x) + g(x)(2|y|)^{1/2} N\}] \, dy > 0,$$

$$\sigma^2 = \int_{\mathcal{I}} dx \int_{-\infty}^{\infty} \int_{-\infty}^{\infty} \int_{-\infty}^{\infty} \mathrm{cov}(\min\{0, y_1^2 H''(x) + g(x) W(y_1)\},$$

$$\min[0, y_2^2 H''(x)$$

$$+ g(x)\{W(y_2 + y_3) - W(y_3)\}]) \, dy_1 \, dy_2 \, dy_3,$$

where the random variable $N$ has the standard normal distribution and $W$ denotes a standard Brownian motion. It is clear that $\mu$ is finite; our proof of Theorem 3.1 will show that $\sigma^2$ is also well defined and finite.

THEOREM 3.1. *Assume the distribution function $F$ has three continuous derivatives on an open interval $\mathfrak{I}'$ which contains the compact bounded $\mathfrak{I}$, and that the density $f = F' > 0$ on $\mathfrak{I}$. If $H \in H_{01}$, then $T = n^{-1}\mu + n^{-7/6}\sigma N_n$, where $N_n$ is asymptotically normally distributed with zero mean and unit variance.*

A version of the theorem continues to hold if the distribution function $F = F_n$ is allowed to depend on $n$. The main requirements in this case are that the regularity conditions hold in a contiguous way, and $F_n$ converge sufficiently fast to a proper limiting distribution, $G$ say. In particular, $F_n$ and $G$ (the former for all sufficiently large $n$) should satisfy the conditions of the theorem, and, for $j = 0$, 1 and 2, $F_n^{(j)} - G^{(j)}$ should converge to 0 at a faster rate than $n^{-1/6}$, uniformly on $\mathfrak{I}'$. Under these assumptions, the limiting distribution of $T$ is that defined when, in the definitions of $\mu$ and $\sigma$, $F$ is replaced by $G$. The proof requires only minor modifications.

This result may be used to prove that if $H \in H_{01}$, and under mild conditions on $h$ and $K$, the bootstrap estimator of the distribution of $T$ is strongly consistent for the limiting distribution of $T$. Our next theorem will state this result. To formulate it, put $\mathcal{H}(\xi_1, \xi_2) = [n^{-\xi_1}, n^{-\xi_2}]$, where

(3.2) $$\tfrac{1}{12} < \xi_2 < \xi_1 < \tfrac{2}{9}.$$

Assume that

(3.3) $K$ is a symmetric, compactly supported probability density with a Hölder-continuous derivative.

Note particularly that bandwidths of size $n^{-1/5}$ are in $\mathcal{H}(\xi_1, \xi_2)$ if (3.2) holds. Indeed, conventional bandwidth selectors, for example, those based on bootstrap methods, cross-validation or plug-in rules, satisfy

(3.4) $$P(C_1 n^{-1/5} < \hat{h} < C_2 n^{-1/5}) \to 1$$
$$\text{as } n \to \infty, \text{ for some } 0 < C_1 < C_2 < \infty.$$

Let $T^*$ denote the version of $T$, defined at (2.2), but with $r = 1$ and $w \equiv 1$, and computed from a sample drawn by sampling randomly from the distribution $\widetilde{F}$ conditional on $\mathcal{X}$. Let $\mu$ and $\sigma$ be as in Theorem 3.1, and write $\Phi$ for the standard normal distribution function.

THEOREM 3.2. *Assume that the (possibly random) bandwidth $h$ lies in $\mathcal{H}(\xi_1, \xi_2)$, where $\xi_1$ and $\xi_2$ satisfy (3.2), and that $K$ satisfies (3.3). Suppose* too that $F$ has four bounded derivatives on an open interval $\mathfrak{I}'$ which contains the compact interval $\mathfrak{I}$, that $f > 0$ on $\mathfrak{I}$ and that $H \in H_{01}$. Then, uniformly in $x$ and with probability 1,

$$P\{n^{7/6}(T^* - n^{-1}\mu)/\sigma \leq x | \mathcal{X}\} \to \Phi(x) \tag{3.5}$$

as $n \to \infty$.

Since the conclusion of Theorem 3.1 may be stated equivalently as

$$P\{n^{7/6}(T - n^{-1}\mu)/\sigma \leq x\} \to \Phi(x),$$

then (3.5) may be interpreted as implying that the bootstrap distribution of $T^*$ converges to the limiting distribution of $T$, provided $H \in H_{01}$.

It should be mentioned too that if a starting bandwidth $\hat{h}$ is chosen using a standard method such as the bootstrap, cross-validation or plug-in, and if the method suggested in Section 2.2 is employed to calculate the critical bandwidth $\hat{h}_{\text{crit}}$, then, under the conditions imposed on $F$ and $K$ in Theorem 3.2, it is true with probability 1 that $\hat{h} = \hat{h}_{\text{crit}}$ for all sufficiently large $n$. That is to say, the iterative process used to define $\hat{h}_{\text{crit}}$ stops at the very first step. This is a consequence of two properties: (i) if $H \in H_{01}$, then $H''$ must, in fact, be bounded above zero on the compact interval $\mathfrak{I}$; and (ii) if a bandwidth of conventional size is used, then $\widetilde{H}''$ converges uniformly to $H''$ on $\mathfrak{I}$ with probability 1. Together (i) and (ii) imply that with probability 1 $\widetilde{H}''$ is bounded above zero for all sufficiently large $n$, and, hence, that $\hat{h} = \hat{h}_{\text{crit}}$ for all sufficiently large $n$.

Furthermore, with probability 1 $\hat{h} \in \mathcal{H}(\xi_1, \xi_2)$ for all sufficiently large $n$. Therefore, when $H \in H_{01}$ the calibration step in Section 2.2 degenerates in asymptotic terms to simply using the standard bandwidth selector, in which case its properties are covered by Theorem 3.2. In particular, using a standard bandwidth selector leads to consistent estimation of the limiting distribution of $T$ when $H \in H_{01}$.

3.3. *Strict monotonicity at all but a finite number of points.* Let $H_{02}$ be the following subset of the class of cumulative hazard functions satisfying $H_0$:

$$H_{02} = \{H : H'' \text{ has two continuous derivatives on } \mathfrak{I}, \text{ and } H'' > 0 \text{ on } \mathfrak{I},$$
$$\text{except for a finite number of distinct points } x_1, \ldots, x_m \in \mathfrak{I},$$
$$\text{where } H'' \text{ vanishes and } H^{(4)} > 0\}. \tag{3.6}$$

We assume $m \geq 1$. Note that it is not possible for $H''$ to vanish at a point $x$, for $H^{(4)}$ to be strictly negative there, and at the same time for the hazard rate to be strictly increasing on sufficiently small intervals containing $x$.

The case of strict monotonicity at all but a finite number of points may fairly be interpreted as the boundary between cases where $H \in H_{01}$ and those where the hazard rate has decreasing parts in the vicinities of points $x_1, \ldots, x_n$. The assumption that $H''(x_i) = 0$ and $H^{(4)}(x_i) > 0$ implies that the hazard rate has a "shoulder" at $x_i$ and is on the verge of decreasing there. Therefore, testing in this context means attempting to identify alternative hypotheses in difficult cases; compare [7]. It offers the opportunity to assess performance against local alternative hypotheses, an opportunity we shall take up in Section 3.5. The opportunity is virtually absent in the setting of Section 3.2.

Let $Z_1, \ldots, Z_m$ be independent random variables, $Z_i$ having the distribution of

$$(3.7) \quad -\int_{-\infty}^{\infty} \int_{-\infty}^{\infty} \min\{0, (\tfrac{1}{2}x^2 y^2 + \tfrac{1}{12} y^4) H^{(4)}(x_i) + g(x_i) W(x+y)\} \, dx \, dy,$$

where $W$ denotes a standard Brownian motion. For simplicity, we shall assume that

(3.8) $\qquad\qquad\qquad$ no $x_i$ is an endpoint of $\mathcal{I}$.

Theorem 3.3 has an analogue in the contrary case; it involves altering the distribution of $Z_i$ when $x_i$ is an endpoint.

THEOREM 3.3. *Assume $F$ has four continuous derivatives on an open interval which contains the compact interval $\mathcal{I}$, and that $f = F' > 0$ on $\mathcal{I}$. Suppose too that $H \in H_{02}$ for points $x_1, \ldots, x_m$ in the definition of that function class, and that (3.8) holds. Then we may write $T = n^{-6/7} \sum_{1 \leq i \leq m} Z_{ni}$, where the joint distribution of $(Z_{n1}, \ldots, Z_{nm})$ converges to that of $(\bar{Z}_1, \ldots, Z_m)$.*

Again, a version of the theorem holds when $F = F_n$ varies with $n$. However, a direct analogue of Theorem 3.2 does not exist in this setting. Essentially, this is because a bandwidth that is sufficiently large to ensure convergence of $\widetilde{H}^{(4)}$ to $H^{(4)}$, and so capture the role of $H^{(4)}(x_i)$ in the definition of the distribution of $Z_i$, is too large to allow sufficiently fast convergence for capturing other features of the limiting distribution. Thus, in the "boundary" case treated by Theorem 3.2, there is not a direct way, based on the estimator $\widetilde{F}$ and using a bandwidth that is asymptotic to a nonrandom quantity, of calibrating the test so as to capture the exact distribution of $T$.

Details behind this claim will be given in Section 5.4. These difficulties persist even if $\widetilde{F}$ is computed using a high-order kernel.

One way of overcoming these difficulties would be to locally model the behavior of $F$ in the neighborhood of points $x$ where $\widetilde{H}''(x)$ was small, rather than leaving estimation there up to the generic estimator $\widetilde{F}$ and to use the model directly to estimate the distributions of $Z_1, \ldots, Z_n$. This approach is rather cumbersome, however, and so, for simplicity we shall not consider it further. Moreover, the problems are largely overcome by the calibration method proposed in Section 2.2, the theory of which we treat next.

3.4. *Calibration based on $\hat{h}_{\mathrm{crit}}$*. The calibration method suggested in Section 2.2 produces a test for which the rejection probability, for $H \in H_{02}$, converges to a number that lies strictly between 0 and 1, and so suffers less from the difficulties noted above. First we describe limiting behavior of the critical bandwidth, $\hat{h}_{\mathrm{crit}}$, in the case $H \in H_{02}$. For simplicity we assume there is only a single point, $x_1$, at which $H''$ vanishes.

Define $c = \frac{1}{2}\int u^2 K(u)\, du$ and

$$
\begin{aligned}
S(q,x,y) &= (cq^2 + \tfrac{1}{2}x^2 + \tfrac{1}{12}y^2)H^{(4)}(x_1) \\
&\quad + q^{-2}g(x_1)\int_{-\infty}^{\infty} K''(u)\, du \\
&\quad \times \int_0^1 \{W(x+ty-qu) + W(x-ty-qu)\}(1-t)\, dt,
\end{aligned}
\tag{3.9}
$$

where $g = f^{1/2}/(1-F)$ and $W$ is a standard Brownian motion. Let $Q > 0$ denote the infimum of values $q > 0$ such that $S(q,x,y) \geq 0$ for all real $x, y$.

THEOREM 3.4. *Assume the conditions of Theorem 3.3, but with $m = 1$. Suppose too that $K$ is a symmetric, compactly supported probability density with two Hölder-continuous derivatives, and that the starting bandwidth $\hat{h}$ used to initiate the algorithm that produces $\hat{h}_{\mathrm{crit}}$ satisfies (3.4). Then $n^{1/7}\hat{h}_{\mathrm{crit}} \to Q$ in distribution as $n \to \infty$.*

Next we describe the asymptotic rejection probability for the test when $H \in H_{02}$. For $0 < \alpha < 1$, define $z_\alpha$ to be the $\alpha$-level quantile of the distribution defined at (3.7) in the case $i = 1$. Noting (3.7), we see that we may write $z_\alpha$ as a continuous function of $H^{(4)}(x_1)$ and $g(x_1)$, say $z_\alpha = \Gamma_\alpha\{H^{(4)}(x_1), g(x_1)\}$. Put

$$
\begin{aligned}
S(x) &= S(Q,x,0) \\
&= (cQ^2 + \tfrac{1}{2}x^2)H^{(4)}(x_1) + Q^{-2}g(x_1)\int_{-\infty}^{\infty} K''(u)W(x-Qu)\, du.
\end{aligned}
\tag{3.10}
$$

It follows from the definition of $Q$ that, with probability 1, (a) $S(x) \geq 0$ for $-\infty < x < \infty$, (b) there exists a unique (random) point $x = A$ at which $S(x) = 0$, and (c) $S'(A) = 0$ and $S''(A) > 0$. [To appreciate why, observe that $S$ is asymptotically proportional to $\widetilde{H}''(x_1 + n^{-1/7}x)$, after taking the bandwidth to equal $n^{-1/7}Q$. Note that the second derivative of $S$ is well defined and continuous as long as $K$ has three continuous derivatives.]

Let $Z_1$ denote the random variable at (3.7) when $i = 1$, constructed using the same Brownian motion $W$ as at (3.10). Therefore, $Z_1$ and $S''(A)$ are linked through $W$. In interpreting the theorem below, note that the probability that $Z_1 \leq \Gamma_\alpha\{H^{(4)}(x_1), g(x_1)\}$ equals $\alpha$.

THEOREM 3.5. *Assume the conditions of Theorem 3.4, but with the additional requirement that $K$ have three continuous derivatives. Take $h = \hat{h}_{\mathrm{crit}}$. Then the rejection probability for the bootstrap test converges as $n \to \infty$ to the probability that $Z_1 \leq \Gamma_\alpha\{S''(A), g(x_1)\}$.*

3.5. *Power against local alternatives and optimality.* Let $F$ denote a four-times continuously-differentiable distribution function for which the corresponding hazard rate is in $H_{02}$. Assume for simplicity that there is only one point at which, for this $F$, $H''$ vanishes on $\mathfrak{I}$. Let this point be $x_1 = 0$, and take it to be an interior point of $\mathfrak{I}$. Since $H \in H_{02}$, then $H^{(3)}(0) = 0$ and $H^{(4)}(0) > 0$.

We shall add a "wiggle" to $F$ in the vicinity of the origin, such that the perturbed distribution violates the null hypothesis. The perturbation will be chosen so that it is only barely detectable using an optimal parametric method, that is, the likelihood-ratio test. We shall then explore the performance of our nonparametric test, based on the statistic $T$, and show that it too is able to detect the wiggle.

The perturbation, $a\varepsilon^4\Psi(x/\varepsilon)$, is based on a four-times continuously-differentiable function $\Psi$ supported on $[-1, 1]$. The constant $a > 0$ represents the height of the wiggle, and $\varepsilon = \varepsilon(n) \to 0$ indicates the extent of the perturbation away from its center, at the origin. We shall choose $\varepsilon$ so small that the perturbation is only barely detectable by the likelihood-ratio test. Our construction of the perturbation ensures that, like the distribution $F$ to which it is added, it has four bounded derivatives near the origin.

The perturbed distribution is

$$F_n(x) = F(x) + a\varepsilon^4\Psi(x/\varepsilon). \tag{3.11}$$

(It is possible, for small $n$, that $F_n$ will be decreasing in some region, but for the choice $\varepsilon = n^{-1/7}$ that we shall make, and under the other regularity conditions of Theorem 3.6, $F_n$ will be nondecreasing on $\mathfrak{I}$ for all sufficiently large $n$.) Let $H_n$ denote the cumulative hazard rate corresponding to $F_n$. If we choose $\Psi$ so that $\Psi(x) \equiv -x^4$ in a neighborhood of the origin, then, for each $a > 0$ and all sufficiently large $n$, $H'_n$ is strictly monotone decreasing in a neighborhood of 0. [This neighborhood is of width $O(\varepsilon)$.] Therefore, $F_n$ fails to satisfy the null hypothesis of an increasing hazard rate.

The density $f_n = F'_n$ satisfies $f_n(x) = f(x) + a\varepsilon^3\psi(x/\varepsilon)$, where $\psi = \Psi'$. Since $f_n$ must be a density, then $\int \psi = 0$. Now,

$$\log\{f_n(x)/f(x)\} = \frac{a\varepsilon^3\psi(x/\varepsilon)}{f(x)} - \frac{a^2\varepsilon^6\psi(x/\varepsilon)^2}{2f(x)^2} + O(\varepsilon^9).$$

Therefore, putting $b(a) = \frac{1}{2}a^2 f(0)^{-1}\int \psi^2$, $f_+ = f_n$ and $f_- = f$, we have, taking the $\pm$ signs, respectively,

$$\int f_\pm(x)\log\{f_n(x)/f(x)\}\,dx = \pm b(a)\varepsilon^7 + o(\varepsilon^7). \tag{3.12}$$

It follows from (3.12) that the expected log-likelihood ratio, for a sample of size $n$, is of size $n\varepsilon^7$. Choosing $\varepsilon$ such that this quantity is bounded away from zero and infinity, in particular, $\varepsilon = n^{-1/7}$, makes the perturbation only barely detectable. In that case, a likelihood-ratio test for discriminating between $f$ and $f_n$ does not have asymptotically perfect accuracy.

Our test is able to detect local alternatives such as $F_n$, provided the function $\pi_{2\alpha}(a)$ for our test satisfies

$$\lim_{a \to \infty} \pi_{2\alpha}(a) = 1. \tag{3.13}$$

If (3.13) holds, then our test shares the optimal performance of the likelihood-ratio test.

To establish (3.13), note first that, for $j = 0, \ldots, 4$,

$$H_n^{(j)} = H^{(j)} + \frac{a\varepsilon^{4-j}\Psi^{(j)}(x/\varepsilon)}{1-F(x)} + O\{\varepsilon^{5-j}I(|x| \leq \varepsilon)\},$$

uniformly in $x$. In particular, the second derivative of $H_n - H$ is of size $n^{-2/7}$, and the fourth derivative is asymptotic to $a\Psi^{(j)}(x/\varepsilon)/\{1-F(x)\}$. Using these properties, and the arguments leading to Theorem 3.5, the following result may be proved. It verifies (3.13) in the case where the test in question is the bootstrap-calibrated one proposed in Section 2.

THEOREM 3.6. *Assume the hazard rate of the four-times continuously-differentiable distribution $F$ lies in $H_{02}$, with $m = 1$ and $x_1 = 0$; and that $F_n$ is given by (3.11), where the four-times continuously-differentiable function $\Psi$ is supported on $[-1, 1]$ and satisfies $\Psi(x) = -x^4$ in a neighborhood of the origin. Suppose too that $\varepsilon$ in (3.11) is $n^{-1/7}$, that $h = \hat{h}_{\mathrm{crit}}$, and that the starting bandwidth $\hat{h}$ satisfies (3.4). Let $p_\alpha(a, n)$ denote the probability that the bootstrap-calibrated test of the null hypothesis of monotone hazard rate rejects the null hypothesis when applied to data from $F_n$. Then (a) $p_\alpha(a, n)$ converges to a limit, $\pi_{2\alpha}(a)$ say, as $n \to \infty$, and (b) $\pi_{2\alpha}(a)$ satisfies (3.13) as $a \to \infty$.*

3.6. *Rejection probability under the null hypothesis, and power against fixed alternatives.* The result below shows that the bootstrap-calibrated form of our test is asymptotically consistent in rejecting the null hypothesis whenever it is violated by a fixed alternative.

THEOREM 3.7. *Assume $F$ has two continuous derivatives on an open interval $\mathfrak{I}'$ which contains the compact interval $\mathfrak{I}$, that $f > 0$ on $\mathfrak{I}$, but that the hazard rate for $F$ is strictly decreasing in a subinterval of $\mathfrak{I}$. Suppose too that $K$ satisfies (3.3), that $K(0) \neq 0$, that $E|X| < \infty$ and that the starting bandwidth $\hat{h}$ for the algorithm leading to $\hat{h}_{\mathrm{crit}}$ defined in Section 2.2 satisfies (3.4). Then $P\{T \geq \hat{c}(\alpha)\} \to 1$, as $n \to \infty$, for each $0 < \alpha < 1$, where $\hat{c}(\alpha)$ is the bootstrap critical point introduced in Section 2.2.*

3.7. *Calibration against the exponential distribution.* Put $A(x) = B\{F(x)\}/\{1-F(x)\}$, where $B$ is a standard Brownian bridge, and define

$$T_0 = \iint_{x,y \colon x+y, x-y \in \mathcal{I}} \max\{0, A(x+y) + A(x-y) - 2A(x)\} \, dx \, dy.$$

In this notation, and using standard Gaussian approximations to the empirical distribution $\widehat{F}$ (see, e.g., [17]), it can be proved that if $F$ is taken to be exponential over $\mathcal{I}$, then $n^{1/2}T \to T_0$ in distribution. This result follows from the fact that, in the exponential case, the cumulative hazard rate is linear. In particular, in that setting $H$ is in neither $H_{01}$ nor $H_{02}$.

Therefore, if we calibrate $T$ by reference to an exponential distribution, then the critical points for the test will be distant $n^{-1/2}$ from the origin. However, if $H$ is in either $H_{01}$ or $H_{02}$, this is much further from zero than the actual critical points of the distribution of $T$. Indeed, we know from Theorems 3.1 and 3.3 that those points are distant only $O(n^{-1})$ from zero when $H \in H_{01}$, and only distant $O(n^{-6/7})$ when $H \in H_{02}$. (The same is true of the bootstrap critical-point estimator suggested in Section 2.2.) It follows that, for each value of the nominal rejection probability of an exponentially calibrated test, the exact rejection probability (for $H$ in either $H_{01}$ or $H_{02}$) will converge to 0 as $n \to \infty$.

Put another way, the exponentially calibrated test will become ultra-conservative as sample size increases. In particular, it will fail, asymptotically, to detect the perturbation-type null hypothesis discussed in Section 3.5. In order for detection to be even barely possible in that setting, the perturbation $\varepsilon^4 \Psi(x/\varepsilon)$ (with $\varepsilon = n^{-1/7}$) would have to be increased by the factor $n^{3/14}$.

**4. Simulations.** Simulations are carried out for two models. First, consider a variable $X$ with hazard rate

(4.1) $$H'(x) = a\{(x-b)^3 + b^3\} + c + dx^2,$$

where $x, a, b, c > 0$ and $d$ is chosen such that $H'(x) > 0$ for all $x > 0$. The distribution function corresponding to this hazard function is given by

$$F(x) = 1 - \exp[-a\{\tfrac{1}{4}(x-b)^4 + b^3 x\} - cx - \tfrac{1}{3}dx^3].$$

It is readily verified that $H \in H_{01}$ when $d > 0$, $H \in H_{02}$ when $d = 0$ and $H$ is in neither $H_{01}$ nor $H_{02}$ when $d < 0$. Figure 1 shows the graph of $H'(x)$ for certain values of the parameters.

The simulations are based on 2000 samples of size $n = 50$ and, for each simulated sample, 2000 resamples are generated. The interval $\mathcal{I}$ on which the test statistic $T$ is based is given by $[0, F^{-1}(0.95)]$. The starting bandwidth $\hat{h}$ is determined from the normal reference rule for plug-in estimation, that is, $\hat{h} = 1.06 n^{-1/5} \hat{\sigma}$, where $\hat{\sigma}$ is the estimated standard error of $X$. The kernel function used is the normal kernel. The results for $a = 2.5, b = 0.75, c = 0.50$, for several values of $d$ and for $\alpha = 0.10$ are presented in Figure 2. The power curve starts at $-1.14$, which is the smallest possible value of $d$ for this choice of parameters. The results for other choices of the parameters and for $\alpha = 0.05$ are similar. For most choices slightly conservative rejection probabilities are observed. As a comparison we also implemented the global sign test of Proschan and Pyke [19] and the local sign test of Gijbels and Heckman [13]. From Figure 2 it is clear that the power curves of both tests are considerably below the curve of the new test. The power of the global test is even identical to zero for all values of $d$. This confirms what was explained in Sections 1 and 3.7 about the lack of power of tests based on calibration with respect to the exponential distribution.

Next, we consider hazard rates which contain a small "bump" and we study how well the three tests are able to detect this little perturbation from $H_0$. The hazard rate considered is

(4.2)  $H'(x) = \exp[\gamma \log x + \beta (2\pi\sigma^2)^{-1/2} \exp\{-(x-\mu)^2/(2\sigma^2)\}],$

where $x, \sigma, \mu > 0$ and $\gamma$ and $\beta$ are real numbers. This model is also considered in [13]. Graphs of this hazard rate for different values of the parameters are presented in Figure 3. It is clear that, for $\beta$ sufficiently large, the hazard rate contains a "bump" at $x = \mu$. The simulation results are obtained from 1000 samples of size $n = 50$ and for the bootstrap procedure 1000 resamples are used. The results are shown in Table 1. Clearly, the hypothesis $H_0$ is only satisfied when $\gamma = 0, 0.50$ or $1$ and $\beta = 0$. In comparison with the local sign test of Gijbels and Heckman [13] and the global sign test of Proschan and Pyke [19], the new testing procedure is now leading to rejection probabilities that are most of the time higher, but not always. Also note that the new test tends to be anticonservative, while the global and local test are, on the contrary, quite conservative. This has to be taken into account when comparing the powers of the three curves.

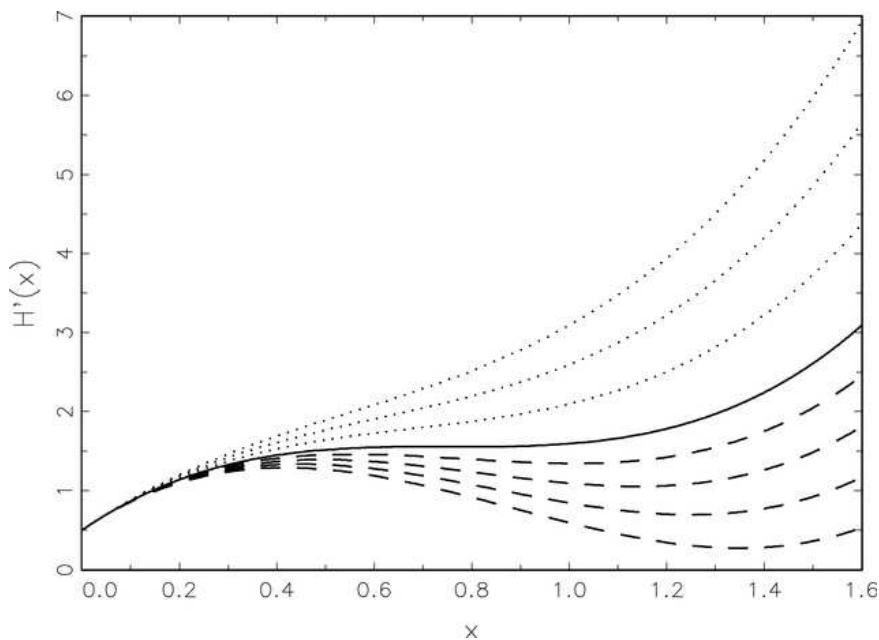

FIG. 1. *Graph of $H'(x)$ for model* (4.1) *when* $a = 2.5, b = 0.75, c = 0.5$ *and* $d = -1, -0.75, -0.50$ *and* $-0.25$ (*dashed curves*), $d = 0$ (*full curve*) *and* $d = 0.5, 1$ *and* $1.5$ (*dotted curves*).

## 5. Technical arguments.

5.1. *Proof of Theorem* 3.1. Define $\Delta_{0F} = \widehat{F} - F$, and observe that

$$(5.1) \quad \widehat{H} = H + \frac{\Delta_{0F}}{1-F} + O_p(n^{-1}), \qquad \Delta_{0F} = n^{-1/2}B(F) + O_p(n^{-1}\log n),$$

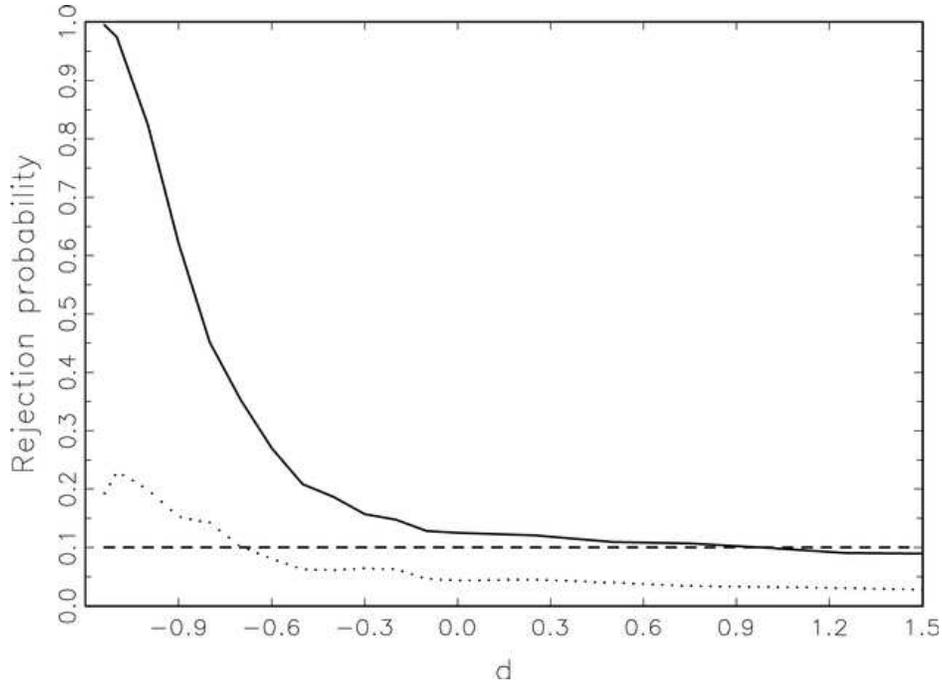

FIG. 2. *Rejection probability for model* (4.1), *when* $a = 2.5, b = 0.75, c = 0.5$, *and for a range of values for* $d$. *The full curve is obtained with the new test, the dotted curve with the local test of Gijbels and Heckman* [13], *while the dashed curve represents the nominal level* $\alpha = 0.10$. *The global test of Proschan and Pyke* [19] *has everywhere zero power.*

where the first result holds uniformly on $\mathfrak{I}$, the second uniformly on the real line and $B$ denotes a Brownian bridge, the construction of which depends on the data. The first identity at (5.1) follows by simple Taylor expansion, while the second uses results of Komlós, Major and Tusnády [17]. Together the identities imply that

$$(5.2) \qquad \widehat{H} = H + n^{-1/2}\frac{B(F)}{1-F} + O_p(n^{-1}\log n),$$

uniformly on $\mathfrak{I}$.

Assume $H \in H_{01}$, and, given a function $\psi(x)$ defined for $x \in \mathfrak{I}$, put $\psi(x,y) = \psi(x+y) + \psi(x-y) - 2\psi(x)$ whenever $x+y, x-y \in \mathfrak{I}$. Now $H(x,y) = y^2 H''(x+\theta y)$, where $-1 \le \theta = \theta(x,y) \le 1$. Hence, for $H \in H_{01}$,

$$\inf_{x,y:\, x+y,x-y\in\mathfrak{I}} y^{-2} H(x,y) > 0.$$

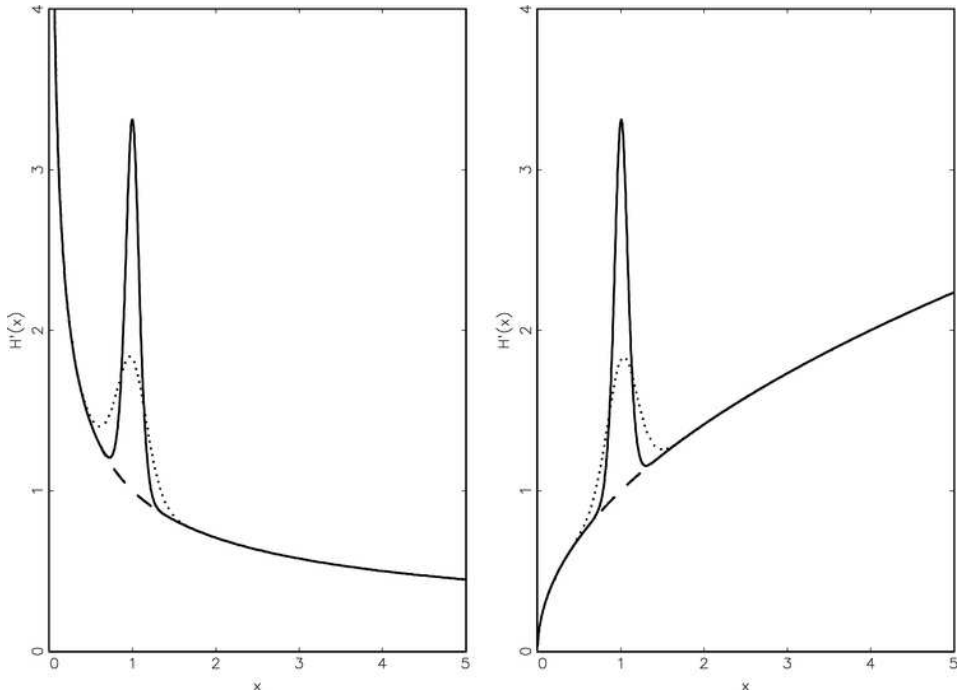

FIG. 3. *Graph of $H'(x)$ for model (4.2) when $\mu = 1$ and $\beta = 0$ (dashed curve), $\beta = 0.3$ and $\sigma = 0.1$ (full curve) and $\beta = 0.3$ and $\sigma = 0.2$ (dotted curve). For the figure on the left $\gamma = -0.5$, for the one on the right $\gamma = 0.5$.*

TABLE 1
*Rejection probability for model* (4.2) *and for* $\alpha = 0.10$. *The numbers in italic are rejection probabilities under the null hypothesis*

| Parameter | Test | \multicolumn{5}{c}{$\gamma$} |
|---|---|---|---|---|---|---|
| | | $-0.50$ | $-0.25$ | $0$ | $0.50$ | $1$ |
| $\beta = 0$ | New | 0.833 | 0.643 | *0.437* | *0.189* | *0.121* |
| | Global | 1.00 | 0.800 | *0.100* | *0.000* | *0.000* |
| | Local | 0.983 | 0.416 | *0.100* | *0.034* | *0.027* |
| $\beta = 0.3$ | New | 0.675 | 0.753 | 0.772 | 0.656 | 0.508 |
| $\sigma = 0.1$ | Global | 0.997 | 0.458 | 0.019 | 0.000 | 0.000 |
| $\mu = 1$ | Local | 0.962 | 0.291 | 0.178 | 0.176 | 0.154 |
| $\beta = 0.3$ | New | 0.715 | 0.714 | 0.663 | 0.443 | 0.277 |
| $\sigma = 0.2$ | Global | 0.999 | 0.588 | 0.035 | 0.000 | 0.000 |
| $\mu = 1$ | Local | 0.968 | 0.301 | 0.114 | 0.065 | 0.054 |

We may therefore deduce from (5.2) that, for some constant $C_1 > 0$, $-\widehat{H}(x,y) > 0$ only if

$$(5.3) \quad C_1 n^{1/2} y^2 \leq \max\{|B\{F(x+y)\} - B\{F(x)\}|, |B\{F(x-y)\} - B\{F(x)\}|\} + n^{-1/2}(\log n) A_n,$$

where the random variable $A_n$ does not depend on $x$ or $y$ and equals $O_p(1)$ as $n \to \infty$.

For each $x$, let $Y(x)$ denote the supremum of values $y$ such that $x+y, x-y \in \mathfrak{I}$ and (5.3) holds. Then for each $x$, $Y(x) = O_p(n^{-1/3})$. Since

$$(5.4) \quad |B(t+u) - B(t)| = O_p(|u \log u|^{1/2})$$

uniformly in $t, u$ such that $0 < t, t+u < 1$, then

$$(5.5) \quad \sup_{\mathfrak{I}} Y(x) = O_p\{(n^{-1} \log n)^{1/3}\}.$$

Defining $\Delta_{1F} = B(F)$ and $\Delta_{2F} = B(F)/(1-F)$, we deduce first by Taylor expansion and then application of (5.4) that

$$(5.6) \quad \begin{aligned} \Delta_{2F}(x,y) &= \frac{B\{F(x+y)\}}{1-F(x)}\left(1 + \frac{yf(x)}{1-F(x)}\right) \\ &\quad + \frac{B\{F(x-y)\}}{1-F(x)}\left(1 - \frac{yf(x)}{1-F(x)}\right) - 2\frac{B\{F(x)\}}{1-F(x)} + O_p\{Y(x)^2\} \\ &= \frac{\Delta_{1F}(x,y)}{1-F(x)} + O_p\{Y(x)^{3/2}(\log n)^{1/2}\}, \end{aligned}$$

uniformly in $x \in \mathfrak{I}$ and $|y| \leq Y(x)$. Therefore,

$$T = \int_{\mathfrak{I}} dx \int_{-Y(x)}^{Y(x)} \max\{0, -\widehat{H}(x,y)\} dy$$

$$
\begin{aligned}
(5.7) \quad &= -\int_{\mathcal{I}} dx \int_{-Y(x)}^{Y(x)} \min\{0, H(x,y) + n^{-1/2}\Delta_{2F}(x,y)\} \, dy \\
&\quad + O_p\{(n^{-1}\log n)^{4/3}\} \\
(5.8) \quad &= -\int_{\mathcal{I}} dx \int_{-Y(x)}^{Y(x)} \min\left\{0, H(x,y) + n^{-1/2}\frac{\Delta_{1F}(x,y)}{1-F(x)}\right\} dy \\
&\quad + O_p\{(n^{-1}\log n)^{4/3}\},
\end{aligned}
$$

where the second identity follows from (5.2) and (5.5), and the third comes from (5.5) and (5.6).

Let $W$ denote the standard Brownian motion through which $B$ may be expressed as $B(t) = W(t) - tW(1)$ for $0 \leq t \leq 1$. Put $\Delta_{3F} = W(F)$. Observe that $\Delta_{1F}(x,y) - \Delta_{3F}(x,y) = O_p\{Y(x)^2\}$ uniformly in $x \in \mathcal{I}$ and $|y| \leq Y(x)$. Therefore, (5.5) and (5.8) imply that

$$
(5.9) \quad T = -\int_{\mathcal{I}} dx \int_{-Y(x)}^{Y(x)} \min\left\{0, H(x,y) + n^{-1/2}\frac{\Delta_{3F}(x,y)}{1-F(x)}\right\} dy \\
+ O_p\{(n^{-1}\log n)^{4/3}\}.
$$

Since $H'''$ is bounded, then $H(x,y) = y^2 H''(x) + O(|y|^3)$ as $y \to 0$, uniformly in $x \in \mathcal{I}$. From this result, (5.5) and (5.9), we deduce that

$$
(5.10) \quad T = T_1 + O_p\{(n^{-1}\log n)^{4/3}\},
$$

where

$$
\begin{aligned}
(5.11) \quad -T_1 &= \int_{\mathcal{I}} dx \int_{-Y(x)}^{Y(x)} \min\left\{0, y^2 H''(x) + n^{-1/2}\frac{\Delta_{3F}(x,y)}{1-F(x)}\right\} dy \\
&= n^{-1} \int_{\mathcal{I}} dx \int_{\mathcal{I}_n(x)} \min\left\{0, z^2 H''(x) + n^{1/6}\frac{\Delta_{3F}(x, n^{-1/3}z)}{1-F(x)}\right\} dz
\end{aligned}
$$

and $\mathcal{I}_n(x)$ denotes the set of $y$ such that both $x + n^{-1/3}y$ and $x - n^{-1/3}y$ lie in $\mathcal{I}$.

Put

$$
W_x(y) = n^{1/6}[W\{F(x) + n^{-1/3}yf(x)\} - W\{F(x)\}]/f(x)^{1/2},
$$

which, like $W$, is a standard Brownian motion. It may be proved from (5.11) that

$$
-nE(T_1) \\
= \int_{\mathcal{I}} dx \int_{-\infty}^{\infty} E\Big\{\min\Big(0, z^2 H''(x) \\
+ \frac{f(x)^{1/2}}{1-F(x)}
$$

$$\times \left[ W_x \left\{ z + \frac{1}{2} n^{-1/3} z^2 f'(x) f(x)^{-1} \right\} \right.$$
$$\left. + W_x \left\{ -z + \frac{1}{2} n^{-1/3} z^2 f'(x) f(x)^{-1} \right\} \right] \right) \right\} dz$$
$$+ o(n^{-1/6}).$$

From this result and the fact that, for $0 < |u| < |z|$, $W_x(z+u) + W_x(-z+u)$ has the normal $N(0, 2|z|)$ distribution, we deduce that

(5.12) $$nE(T_1) = \mu + o(n^{-1/6}),$$

where $\mu$ is as defined in Section 3.

To derive a central limit theorem for $T_1$, we first approximate $T_1$ by a sum of 3-dependent random variables, as follows. Define $\lambda_n = \log n$ and $\delta = \delta(n) = \lambda_n (n^{-1} \log n)^{1/3}$. Put

$$-T_2 = \int_{\mathcal{I}} dx \int_{-\delta}^{\delta} \min \left\{ 0, y^2 H''(x) + n^{-1/2} \frac{\Delta_{3F}(x,y)}{1 - F(x)} \right\} dy,$$

$$-T_2(i) = \int_{\mathcal{I} \cap (i\delta, (i+1)\delta)} dx \int_{-\delta}^{\delta} \min \left\{ 0, y^2 H''(x) + n^{-1/2} \frac{\Delta_{3F}(x,y)}{1 - F(x)} \right\} dy;$$

compare these definitions with the first identity at (5.11). Then $T_2 = \sum_i T_2(i)$. Note that, since Brownian motion has independent increments, $T_2(i)$ is stochastically independent of $T_2(j)$ for $|i - j| \geq 3$.

In view of (5.5), the probability that $\max_{\mathcal{I}} |Y(x)| \leq \delta$ converges to 1 as $n \to \infty$. Note too that $\max_{x \in \mathcal{I}} |Y(x)| \leq \delta$ implies $T_1 = T_2$. Hence, if we prove that the following three results are true: (a) $\text{var } T_2 \sim \text{var } T_1 \sim \sigma^2 n^{-7/3}$, (b) $(T_2 - ET_2)/(\text{var } T_2)^{1/2}$ has an asymptotic standard normal distribution, and (c) $n^{7/6}(ET_1 - ET_2) \to 0$; then it will follow that $n^{7/6}(T_1 - ET_1)/\sigma$ has an asymptotic standard normal distribution. Theorem 3.1 is a consequence of this property and (5.12).

Result (c) may be proved using the argument leading to (5.12), and the first asymptotic relation in (a) may be derived using the method giving the second. Therefore, it suffices to show that (b) holds and that (d) $\text{var } T_1 \sim \sigma^2 n^{-7/3}$.

To prove (b), let $C > 0$ and define $T_3(i) = n^{7/6} \delta^{-1/2} T_2(i)$, $T_4(i) = T_3(i) \times I\{|T_3(i)| \leq C\}$, $T_5(i) = T_3(i) - T_4(i)$ and $T_j = \sum_i T_j(i)$ for $j = 4, 5$. For all sufficiently large $C$, the variance of $T_4$, and the number of nondegenerate summands $T_4(i)$, are both asymptotic to constant multiples of $\delta^{-1}$; and the summands are uniformly bounded. Therefore, using a central limit theorem for uniformly bounded $m$-dependent random variables (see, e.g., Theorem 7.3.1, page 214 of [9]), we may prove that $(T_4 - ET_4)/(\text{var } T_4)^{1/2}$ has an asymptotic standard normal distribution; call this result (e). The argument

that we shall use to prove (d) may be employed to show that as $C \to \infty$, (f) $\lim_{n\to\infty} \delta \times \operatorname{var} T_4 \to \sigma^2$ and (g) $\lim_{n\to\infty} \delta \operatorname{var} T_5 \to 0$. Combining (e)–(g), we deduce that $(T_3 - ET_3)/(\operatorname{var} T_3)^{1/2}$ has an asymptotic normal distribution. This is equivalent to (b).

It remains to derive (d). Recall that $g = f^{1/2}/(1-F)$, and define

$$U_x(y) = n^{1/6}[W\{F(x+n^{-1/3}y)\} - W\{F(x)\}]/f(x)^{1/2},$$
$$V_x(y) = U_x(y) + U_x(-y),$$
$$W_1(x_1, y_1) = \min\{0, y_1^2 H''(x_1) + g(x_1)V_{x_1}(y_1)\},$$
$$W_2(x_1, x, y_2) = \min\{0, y_2^2 H''(x_1 + n^{-1/3}x) + g(x_1 + n^{-1/3}x)V_{x_1 + n^{-1/3}x}(y_2)\}$$

and $\mathfrak{J}_n(x_1) = \{x : x_1 + n^{-1/3}x \in \mathfrak{J}\}$. In this notation,

$$n^2 \operatorname{var} T_1 = \int_{\mathfrak{J}} dx_1 \int_{\mathfrak{J}} dx_2 \int_{\mathfrak{J}_n(x_1)} dy_1$$
$$\times \int_{\mathfrak{J}_n(x_2)} \operatorname{cov}[\min\{0, y_1^2 H''(x_1) + g(x_1)V_{x_1}(y_1)\},$$
$$\min\{0, y_2^2 H''(x_2) + g(x_2)V_{x_2}(y_2)\}] dy_1 dy_2$$
$$= n^{-1/3} \int_{\mathfrak{J}} dx_1 \int_{\mathfrak{J}_n(x_1)} dx \int_{\mathfrak{J}_n(x_1)} dy_1$$
$$\times \int_{\mathfrak{J}_n(x_1+n^{-1/3}x)} \operatorname{cov}\{W_1(x_1,y_1), W_2(x_1,x,y_2)\} dy_2,$$

where $\mathfrak{J}_n(x)$ is as defined below (5.11). In view of the independent increments of Brownian motion, the random variables $V_{x_1}(y_1)$ and $V_{x_1+n^{-1/3}x}(y_2)$ are independent if $|y_1| + |y_2| \leq |x|$. In this case, the covariance in the second identity above vanishes. Therefore,

(5.13)
$$n^{7/3} \operatorname{var} T_1 = \int_{\mathfrak{J}} dx_1 \int_{\mathfrak{J}_n(x_1)} dx$$
$$\times \iint_{y_1,y_2:\,|y_1|+|y_2|>|x|;\mathcal{C}(x_1,x)} \operatorname{cov}\{W_1(x_1,y_1), W_2(x_1,x,y_2)\} dy_1 dy_2,$$

where $\mathcal{C}(x_1, x)$ denotes the constraint that $y_1 \in \mathfrak{J}_n(x_1)$ and $y_2 \in \mathfrak{J}_n(x_1 + n^{-1/3}x)$.

The random variables $|W_1(x_1, y_1)|$ and $|W_2(x_1, x, y_2)|$ are respectively bounded by $C_1|N_1|I(|N_1| > C_2 y_1^2)$ and $C_1|N_2|I(|N_2| > C_2 y_2^2)$, where $N_1$ and $N_2$ are standard normal random variables, and $C_1$ and $C_2$ are positive constants not depending on $x_1$, $x$, $y_1$ or $y_2$, although the correlation between

$N_1$ and $N_2$ does depend on these quantities. We may therefore deduce that, for constants $C_3, C_4 > 0$,

$$
\begin{aligned}
&|\operatorname{cov}\{W_1(x_1,y_1), W_2(x_1,x,y_2)\}| \\
&\quad \leq C_1 E\{|N_1|^2 I(|N_1| > C_2 y_1^2)\}^{1/2} E\{|N_2|^2 I(|N_2| > C_2 y_2^2)\}^{1/2} \\
&\quad \leq C_3 \exp\{-C_4(y_1^4 + y_2^4)\}.
\end{aligned}
\tag{5.14}
$$

Therefore, $|\operatorname{cov}\{W_1(x_1,y_1), W_2(x_1,x,y_2)\}|$ is bounded above by a function which does not depend on $n$ and whose integral over $-\infty < x < \infty$ and over all real $y_1, y_2$ that satisfy $|y_1| + |y_2| > |x|$ is bounded uniformly in $x_1 \in \mathcal{I}$.

Furthermore, if $V$ is a standard Brownian motion, then

$$
\begin{aligned}
\operatorname{cov}&\{W_1(x_1,y_1), W_2(x_1,x,y_2)\} \\
&\to \operatorname{cov}(\min\{0, y_1^2 H''(x_1) + g(x_1)V(y_1)\}, \\
&\qquad\qquad \min[0, y_2^2 H''(x_1) + g(x_1)\{V(x+y_2) - V(x)\}]),
\end{aligned}
$$

uniformly in $x_1 \in \mathcal{I}$ and $x$, $y_1$ and $y_2$ in any compact set. We may therefore deduce from (5.13) and the dominated convergence theorem that $\operatorname{var} T_1 \sim \sigma^2 n^{-7/3}$, which is the desired result (d). Note that (5.14) also implies the finiteness of $\sigma^2$.

5.2. *Proof of Theorem* 3.2. Put $\nu_0 = 0$ and $\nu_j = 2j - 1$ for $j \geq 1$. Observe that, for $j = 0, 1, 2$ and each $\eta > 0$,

$$
\widetilde{F}^{(j)}(x) - F^{(j)}(x) = O_p\{(nh^{\nu_j})^{\eta - (1/2)} + h^2\},
\tag{5.15}
$$

uniformly in $h \in \mathcal{H}(\xi_1, \xi_2)$ and $x \in \mathcal{I}'$. (The assumption that $F$ has four bounded derivatives is needed to derive the $O_p(h^2)$ remainder term in (5.15) when $j = 2$. The other part of the remainder at (5.15), which applies to the error of the left-hand side about its mean, may be obtained by applying the stochastic approximation of Komlós, Major and Tusnády [17].) It follows from this property and (2.3) that, with probability 1, $\widetilde{H}''$ converges to $H''$ uniformly in $h \in \mathcal{H}(\xi_1, \xi_2)$ and $x \in \mathcal{I}'$. We may choose $\mathcal{I}'$ and $\varepsilon > 0$ such that $H'' > \varepsilon$ on $\mathcal{I}'$. In this case, and with probability 1, $\widetilde{H}'' > \frac{1}{2}\varepsilon$ on $\mathcal{I}'$ for all sufficiently large $n$. In particular, for all sufficiently large $n$, the hazard rate corresponding to $\widetilde{F}$ lies in $H_{01}$.

The argument leading to Theorem 3.1 may now be used to prove that (3.5) holds when $\widetilde{F}$, rather than $F$, is the sampled distribution, provided $\mu$ and $\sigma$ at (3.5) are replaced by the analogous functionals of $\widetilde{F}$. Let these be $\widetilde{\mu}$ and $\widetilde{\sigma}$, respectively, and denote by (R) the corresponding version of (3.5). By (5.15),

$$
|\widetilde{\mu} - \mu| + |\widetilde{\sigma} - \sigma| = O_p\{(nh^3)^{\eta - (1/2)} + h^2\} = o_p(n^{-1/6}),
\tag{5.16}
$$

the second identity holding uniformly in $h \in \mathcal{H}(\xi_1, \xi_2)$ and following from (3.2). Property (3.5) follows from (5.16) and (R).

We should mention that the assumption in Theorem 3.1 that $F$ have three derivatives is imposed for simplicity, and is a little more stringent than necessary. At (5.10) we need only two derivatives and a Hölder condition of order $\frac{1}{2} + \varepsilon$ on $F''$, in which case the $O_p$ term at (5.10) becomes $O_p\{(n^{-1}\log n)^{(3+\varepsilon)/3}\} = o_p(n^{-7/6})$ (as required), the identity holding provided $\varepsilon > 0$. An empirical version of this argument can be developed provided $h \in \mathcal{H}(\xi_1, \xi_2)$ and $\xi_1, \xi_2$ satisfy (3.2).

5.3. *Proof of Theorem* 3.3. The assumption that the hazard rate is nondecreasing and that $H''(x_i) = 0$ implies that $H'''(x_i) = 0$ for $1 \leq i \leq m$. To appreciate why, observe that

$$H(x,y) = y^2 H''(x) + \tfrac{1}{12} y^4 H^{(4)}(x + \theta y),$$

where $-1 \leq \theta = \theta(x,y) \leq 1$. Taking $x = x_i + u$, where $|u|$ is small, and Taylor-expanding, we deduce that

$$H(x_i + u, y) = y^2 u H'''(x_i) + (\tfrac{1}{2} u^2 y^2 + \tfrac{1}{12} y^4) H^{(4)}\{x + \theta'(|u| + |y|)\},$$

where $-1 \leq \theta' \leq 1$. If $H'''(x_i) \neq 0$, then, taking $|u| = |y|^{3/2}$ and choosing the sign of $u$ such that $uH'''(x_i) < 0$, we find that as $y \to 0$, $H(x_i + u, y) = -|y|^{7/2}|H'''(x_i)| + o(|y|^{7/2})$. This implies that $H$ is nonconvex near $x_i$, and so contradicts the assumption that the hazard rate is nondecreasing.

Result (5.2) continues to hold in the setting of Theorem 3.3, and so by (5.7),

(5.17) $$T = T_2 + O_p\{(n^{-1}\log n)^{8/7}\},$$

where

$$T_2 = -\int_{\mathcal{I}} dx \int_{-Y(x)}^{Y(x)} \min\{0, y^2 H''(x) + \tfrac{1}{12} y^4 H^{(4)}(x+\theta y) + n^{-1/2} \Delta_{2F}(x,y)\} \, dy$$

and we redefine $Y(x)$ to equal the supremum of values $y$ such that $x + y, x - y \in \mathcal{I}$ and

$$y^2 H''(x) + \tfrac{1}{12} y^4 H^{(4)}(x + \theta y) + n^{-1/2} \Delta_{2F}(x,y) + n^{-1}(\log n) A_n \leq 0,$$

where the random variable $A_n = O_p(1)$ does not depend on $x$ or $y$. In deriving (5.17), we have used the fact that, by employing arguments leading to (5.5), it may be proved that

$$\sup_{x \in \mathcal{I}} Y(x) = O_p\{(n^{-1}\log n)^{1/7}\}.$$

More analogously to (5.5), it may be shown that if $\eta > 0$ and $\mathcal{J} = \mathcal{J}(\eta)$ is the subset of $\mathcal{I}$ all of whose points are distant at least $\eta$ from each $x_i$, then, using the new definition of $Y(x)$,

$$\sup_{x \in \mathcal{J}} Y(x) = O_p\{(n^{-1}\log n)^{1/3}\}.$$

Using this result and the arguments leading to (5.9) and (5.10), we may show that, if $T_2(\mathcal{J})$ denotes the contribution to $T_2$ from the integral over $x \in \mathcal{J}$, rather than $x \in \mathcal{I}$, then $T_2(\mathcal{J}) = T_3(\mathcal{J}) + o_p(n^{-1})$, where

$$T_3(\mathcal{J}) = -\int_{\mathcal{J}} dx \int_{-Y(x)}^{Y(x)} \min\left\{0, y^2 H''(x) + n^{-1/2}\frac{\Delta_{3F}(x,y)}{1-F(x)}\right\} dy.$$

The methods leading to (5.12) give that $E\{T_3(\mathcal{J})\} = O(n^{-1})$. Therefore,

(5.18) $$T_2(\mathcal{J}) = O_p(n^{-1}).$$

Let $\eta > 0$ be less than half the minimum of $x_{i+1} - x_i$ over $0 \leq i \leq m$, where $x_0$ denotes the lower limit of $\mathcal{I}$ and $x_{m+1}$ is the upper limit. Write $T_2(x_i, \eta)$ for the contribution to $T_2$ from the integral over $x_i - \eta < x < x_i + \eta$. Then

$$T_2(x_i, \eta) = -\int_{-\eta}^{\eta} du \int_{-Y(x_i+u)}^{Y(x_i+u)} \min[0, (\tfrac{1}{2}u^2 y^2 + \tfrac{1}{12}y^4)$$
$$\times H^{(4)}\{x_i + \theta_i(|u|+|y|)\} + n^{-1/2}\Delta_{2F}(x_i, y)]\, dy,$$

where $-1 \leq \theta_i \leq 1$. Changing variables from $(u, y)$ to $(v, z)$, where $u = n^{-1/7}v$ and $y = n^{-1/7}z$, we deduce that

(5.19) $$T_2(x_i, \eta) = -n^{-6/7}\int_{-\infty}^{\infty} dv \int_{-\infty}^{\infty} \min\{0, (\tfrac{1}{2}v^2 z^2 + \tfrac{1}{12}z^4)H^{(4)}(x_i) + g(x_i)W_i(v+z)\}\, dz$$
$$+ o_p(n^{-6/7}),$$

where $W_i$ is a standard Brownian motion. The processes $W_i$, $1 \leq i \leq m$, may be taken to be independent without violating (5.19). Theorem 3.3 now follows on combining (5.18) and (5.19).

5.4. *Reasons for failure of bootstrap version of Theorem* 3.3. In order for $\widetilde{H}^{(4)}$ to consistently estimate $H^{(4)}$, it is necessary that the bandwidth $h$ used to construct $\widetilde{F}$ be of larger order than $n^{-1/7}$. For simplicity, we shall assume below that $h$ is at least of size $n^{\xi-(1/7)}$ for some $\xi > 0$, although our argument may by pursued to an unaltered conclusion when the increase of $h$ over $n^{-1/7}$ is by only a logarithmic factor.

Put $c = \tfrac{1}{2}\int u^2 K(u)\, du$. Observe that, for each $\eta > 0$, $\widetilde{F}'' = F'' + ch^2 F^{(4)} + O_p\{(nh^3)^{\eta-(1/2)}\} + o(h^2)$, uniformly in $x \in \mathcal{I}'$. [Here we have used the fact that $h \geq n^{\xi-(1/7)}$.] It follows that $\widetilde{H}'' = D^2 A(F + ch^2 F'') + O_p\{(nh^3)^{\eta-(1/2)}\} + o(h^2)$, uniformly in $x \in \mathcal{I}'$, where $A(u) = -\log(1-u)$ and $D$ is the differential operator. Now, $D^2 A(F + ch^2 F'') = D^2 A(F) + ch^2 D^2\{F'' A'(F)\} + o(h^2)$, and $D^2\{F''A'(F)\} = D^2\{D^2 A(F) - (F')^2 A''(F)\} = D^2\{H'' - (H')^2\}$. Therefore, $D^2\{F''A'(F)\} = H^{(4)} - 2\{H'H''' + (H'')^2\}$. Hence,

(5.20) $$\widetilde{H}'' = H'' + ch^2[H^{(4)} - 2\{H'H''' + (H'')^2\}]$$
$$+ O_p\{(nh^3)^{\eta-(1/2)}\} + o(h^2),$$

uniformly in $x \in \mathcal{I}'$.

The term of order $(nh^3)^{\eta-(1/2)}$ on the right-hand side of (5.20) is, of course, the result of stochastic error, and performance would only improve if it could be dropped. Let us assume this can be done. Then we can estimate $H''(x)$ with error equal to

$$(5.21) \qquad ch^2[H^{(4)}(x) - 2\{H'(x)H'''(x) + H''(x)^2\}] + o(h^2).$$

Now, the limiting distribution of $T$, when $H \in H_{02}$, is determined by properties of $H$ in arbitrarily small neighborhoods of the points $x_i$, and so it is there that we are most interested in properties of $\widetilde{H}''$. If $x$ is in a decreasingly small neighborhood of $x_i$, the expansion at (5.21) equals

$$ch^2[H^{(4)}(x_i) - 2\{H'(x_i)H'''(x_i) + H''(x_i)^2\}] + o(h^2) = ch^2 H^{(4)}(x_i) + o(h^2),$$

the second identity holding since $H''(x_i) = H'''(x_i) = 0$. Therefore, if we ignore stochastic fluctuations (which are asymptotically equally likely to increase or decrease the value of $\widetilde{H}''$), $\widetilde{H}''(x)$ is distant at least order $h^2$ strictly above zero when $x$ is in the neighborhood of $x_i$. Since $h$ is at least of order $n^{-1/7}$, then the distance of $\widetilde{H}''$ above zero, in the neighborhood of $x_i$, is [with probability at least $\frac{1}{2} + o(1)$] no less than a certain fixed constant multiple of $n^{-2/7}$; call this result (R).

Let $\widehat{H}^*$ denote the bootstrap version of $\widehat{H}$, and recall from the proof of Theorem 3.2 that that limit result derives entirely from fluctuations of $\widehat{H}(x,y)$ below zero when $x$ is close to $x_i$ and $y$ is near zero. If $H \in H_{02}$, these fluctuations occur with a probability that is bounded away from zero as $n$ increases. The perturbations of $\widehat{H}^* - \widetilde{H}$ are of order only $n^{-1/2}$, and, in particular, are of strictly smaller order than $n^{-2/7}$. This property and result (R) imply that the probability that the empirical fluctuations of $\widehat{H}^*$ near $x_1, \ldots, x_m$ ever protrude below zero converges to zero as $n \to \infty$. In consequence, the limit results described by Theorem 3.2 do not apply in the bootstrap setting.

5.5. *Proof of Theorem* 3.4. Observe that

$$(5.22) \quad \begin{aligned} &\widetilde{H}(x+y) + \widetilde{H}(x-y) - 2\widetilde{H}(x) \\ &= y^2 \int_0^1 \{\widetilde{H}''(x+ty) + \widetilde{H}''(x-ty)\}(1-t)\,dt. \end{aligned}$$

Let $H_G$ denote the version of $H$ that arises if $F$ is replaced by a distribution $G$, and note that, by (2.3) and approximations based, for example, on the Komlós, Major and Tusnády [17] embedding,

$$(5.23) \qquad \widetilde{H}'' = H''_{E(\widetilde{F})} + (1-F)^{-1}(\tilde{f}' - E\tilde{f}') + O_p\{(nh)^{\eta-(1/2)}\},$$

uniformly in $x \in \mathcal{I}'$ and in $h \in \mathcal{H}$, for each $\eta > 0$. The argument in Section 5.4 [see particularly (5.20)] shows that

$$H''_{E(\widetilde{F})} = H'' + ch^2[H^{(4)} - 2\{H'H''' + (H'')^2\}] + o(h^2)$$

uniformly on $\mathcal{I}'$ and in $h \in \mathcal{H}$. Therefore,

$$
\begin{aligned}
\int_0^1 &\{H''_{E(\widetilde{F})}(x+ty) + H''_{E(\widetilde{F})}(x-ty)\}(1-t)\,dt \\
&= \int_0^1 \{H''(x+ty) + H''(x-ty)\}(1-t)\,dt + ch^2 H^{(4)}(x_1) + o(h^2) \\
&= H''(x) + (\tfrac{1}{12}y^2 + ch^2)H^{(4)}(x_1) + o(h^2 + y^2)
\end{aligned}
\tag{5.24}
$$

uniformly in $h \in \mathcal{H}$, $|x - x_1| \leq \delta(n)$ and $|y| \leq \delta(n)$ for any sequence $\delta(n) \downarrow 0$.

Furthermore,

$$
\begin{aligned}
\tilde{f}'(x) - E\tilde{f}'(x) &= h^{-2}\int K''(u)\{\widehat{F}(x-hu) - F(x-hu)\}\,du \\
&= h^{-2}n^{-1/2}\int K''(u)[W\{F(x-hu)\} - W\{F(x)\}]\,du \\
&\quad + O_p\{h^{-1}n^{-1/2}(\log n)^{1/2}\},
\end{aligned}
$$

uniformly in $h \in \mathcal{H}$ and $x \in \mathcal{I}'$, where $W$ is a standard Brownian motion. Put $h = n^{-1/7}q$, $x = x_1 + n^{-1/7}s + ty$ and $y = n^{-1/7}z$, and recall that $g = f^{1/2}/(1 - F)$. Then there exists a standard Brownian motion $V$ such that

$$
\begin{aligned}
h^{-2}&n^{-1/2}\{1 - F(x)\}^{-1}[W\{F(x-hu)\} - W\{F(x)\}] \\
&= n^{-2/7}q^{-2}g(x_1)\{V(s+tz-qu) - V(s+tz)\} \\
&\quad + O_p\{n^{-5/14}(\log n)^{1/2}(|q| + |s| + |z|)\},
\end{aligned}
$$

uniformly in $0 \leq t \leq 1$, $|u| \leq C$ for any $C > 0$ and $q, s, z$ such that $n^{-1/7}q \in \mathcal{H}$, $|s| \leq n^{1/7}\delta(n)$ and $|z| \leq n^{1/7}\delta(n)$. Therefore, defining $M = (1-F)^{-1}(\tilde{f}' - E\tilde{f}')$, we have

$$
\begin{aligned}
\int_0^1 &\{M(x+ty) + M(x-ty)\}(1-t)\,dt \\
&= n^{-2/7}q^{-2}\,g(x_1)\int K''(u)\,du \\
&\quad \times \int_0^1 \{V(s+tz-qu) + V(s-tz-qu)\}(1-t)\,dt \\
&\quad + O_p\{n^{-5/14}(\log n)^{1/2}(|q| + |q|^{-1} + |s| + |z|)\},
\end{aligned}
\tag{5.25}
$$

uniformly in $n^{-1/7}q \in \mathcal{H}$, $|s| \leq n^{1/7}\delta(n)$ and $|z| \leq n^{1/7}\delta(n)$.

Combining (5.22)–(5.25), and taking $x = x_1 + n^{-1/7}s$ and $y = n^{-1/7}z$, we deduce that

$$
\begin{aligned}
n^{4/7}\{&\widetilde{H}(x+y) + \widetilde{H}(x-y) - 2\widetilde{H}(x)\} \\
&= z^2(cq^2 + \tfrac{1}{2}s^2 + \tfrac{1}{12}z^2)H^{(4)}(x_1) \\
&\quad + z^2 q^{-2} g(x_1) \int K''(u)\,du \\
&\quad \times \int_0^1 \{V(s+tz-qu) + V(s-tz-qu)\}(1-t)\,dt \\
&\quad + O_p\{n^{-1/14}(\log n)^{1/2} z^2(|q| + |q|^{-1} + |s| + |z|)\} \\
&\quad + o_p\{z^2(q^2 + s^2 + z^2)\},
\end{aligned}
\tag{5.26}
$$

uniformly in $n^{-1/7}q \in \mathcal{H}$, $|s| \leq n^{1/7}\delta(n)$ and $|z| \leq n^{1/7}\delta(n)$. The theorem follows from (5.26).

5.6. *Proof of Theorem* 3.5. Dividing both sides of (5.26) by $z^2$ and letting $z \to 0$, we deduce that, when $h = \hat{h}_{\text{crit}}$, $n^{2/7}\widetilde{H}''(x_1 + n^{-1/7}s) = S(s) + o_p(1)$, uniformly in $|s| \leq n^{1/7}\delta(n)$, where $S(s)$ is defined as at (3.10). Thus, the bootstrap calibration step involves sampling from a distribution whose cumulative hazard rate $\bar{H}$ is convex on $\mathcal{I}$ and satisfies $\bar{H}''(x) > 0$ for all $x \in \mathcal{I}$, excepting a single point $x$ which may be expressed as $x = x_1 + n^{-1/7}A + o_p(n^{-1/7})$, where $A$ is uniquely defined by $S(A) = 0$. At this point $\bar{H}''$ vanishes. Reworking the proof of Theorem 3.3, we deduce that the critical point $\hat{c}(\alpha)$ of $T^*$, defined conditional on the data $\mathcal{X}$, equals $n^{-6/7}\Gamma_\alpha\{S''(A), g(x_1)\} + o_p(n^{-6/7})$. [Here, $T^*$ denotes the value of $T$ computed from an $n$-sample drawn from the distribution $\widetilde{F}(\cdot|\hat{h}_{\text{crit}})$.]

The distribution of $T$ may be represented, in asymptotic form, as before, and in terms of the same Brownian motion $W$ that was used to construct the representation for $\widetilde{H}$ at (5.26). In particular, the Brownian motion $W_1$ appearing at (5.19) (when $i=1$) may be taken identical to the process $V$ at (5.26). Letting $W$ denote the common process, we see that the inequality $T \leq \hat{c}(\alpha)$ may equivalently be written as

$$n^{-6/7}Z_1 + o_p(n^{-6/7}) \leq n^{-6/7}\Gamma_\alpha\{S''(A), g(x_1)\} + o_p(n^{-6/7}), \tag{5.27}$$

where $Z_1$ is defined by (3.7) with $i = 1$. Theorem 3.5 follows from (5.27).

5.7. *Proof of Theorem* 3.7. If the hazard rate $H'$ is not increasing on $\mathcal{I}$, then, for some $\varepsilon > 0$, there exists a nondegenerate rectangle $\mathcal{R}$ such that for all $(x,y) \in \mathcal{R}$, both $x+y$ and $x-y$ lie in $\mathcal{I}$ and $H(x,y) \leq -\varepsilon$. Under the hypotheses of the theorem, $\widehat{H}(x,y) = H(x,y) + o_p(1)$ uniformly in $(x,y) \in \mathcal{R}$, and so $T \geq \varepsilon|\mathcal{R}| + o_p(1)$, where $|\mathcal{R}|$ denotes the area of $\mathcal{R}$. Therefore, the theorem will follow if we prove that, for each $\alpha \in (0,1)$, the point $\hat{c}(\alpha)$ derived using the bootstrap argument in Section 2.2 satisfies

$$P\{\hat{c}(\alpha) > \eta\} \to 0 \quad \text{for each } \eta > 0. \tag{5.28}$$

As $h \to \infty$, $E\{\tilde{f}(x|h)\} = h^{-1}K(0) + o(h^{-1})$ and $E|\tilde{f}'(x|h)| = h^{-3}K''(0) \times E|x - X| + o(h^{-3})$, uniformly in $x \in \mathcal{I}$. Hence, there exists $h_0 > 0$ such that $\{E\tilde{f}(x|h_0)\}^2 \geq 2E|\tilde{f}'(x|h_0)|$ for all $x \in \mathcal{I}$. It may be proved from this property that, with probability converging to 1, $\tilde{f}(x|h_0)^2 \geq |\tilde{f}'(x|h_0)|$ for all $x \in \mathcal{I}$. Therefore, by (2.3), the probability that $\tilde{H}''(x|h_0) \geq 0$ for all $x \in \mathcal{I}$ converges to 1 as $n \to \infty$, and so

(5.29) $$P(\hat{h}_{\text{crit}} \leq h_0) \to 1.$$

Standard calculations of the expected value of a kernel distribution estimator show that, under the conditions of the theorem, for each $h_1 > 0$, there exists $\varepsilon(h_1) > 0$ such that, for all sufficiently large $n$,

$$1 - E\{\widetilde{F}(x|h)\} \geq \varepsilon(h_1) \quad \text{for all } x \in \mathcal{I} \text{ and all } h \in (0, h_1].$$

By employing a stochastic approximation based on the results of Komlós, Major and Tusnády [17], it may be proved that, for each $h_1 > 0$,

$$|\widetilde{F}(x|h) - E\{\widetilde{F}(x|h)\}| = o_p(1) \quad \text{uniformly in } x \in \mathcal{I} \text{ and } h \in (0, h_1].$$

Therefore,

$$P\{1 - \widetilde{F}(x|h) \geq \tfrac{1}{2}\varepsilon(h_1) \text{ for all } x \in \mathcal{I} \text{ and all } h \in (0, h_1]\} \to 1.$$

This result and (5.29) imply that

(5.30) $$P\{1 - \widetilde{F}(x|\hat{h}_{\text{crit}}) \geq \tfrac{1}{2}\varepsilon(h_0) \text{ for all } x \in \mathcal{I}\} \to 1.$$

If $\widehat{F}_h^*$ denotes the bootstrap version of $\widehat{F}$, computed from an $n$-sample drawn from the distribution $\widetilde{F}(\cdot|h)$ rather than from $F$, then for all $\lambda > 0$,

$$\sup_{x \in \mathcal{I}} \sup_{h \in (0, h_1]} E\{|\widehat{F}_h^*(x) - \widetilde{F}(x|h)|^\lambda\} = O(n^{-\lambda/2})$$

for all $\lambda > 0$. (The method of proof involves only direct calculation of moments, first conditional on the data and then unconditionally.) Therefore, if $\mathcal{A}_1$ and $\mathcal{A}_2$ are subsets of $\mathcal{I}$ and $[Cn^{-1/5}, h_1]$, respectively, each of which contains no more than $O(n^D)$ elements, then for each $\lambda > 0$ and by Hölder's inequality,

$$E\left\{\sup_{x \in \mathcal{A}_1} \sup_{h \in \mathcal{A}_2} |\widehat{F}_h^*(x) - \widetilde{F}(x|h)|\right\}$$
$$\leq \left[\sum_{x \in \mathcal{A}_1} \sum_{h \in \mathcal{A}_2} E\{|\widehat{F}_h^*(x) - \widetilde{F}(x|h)|^\lambda\}\right]^{1/\lambda}$$
$$= \{O(n^{2D-(\lambda/2)})\}^{1/\lambda}$$
$$= O(n^{(2D/\lambda)-(1/2)}).$$

Since $D/\lambda$ can be made arbitrarily small by choosing $\lambda$ sufficiently large, then we have proved that, for each $\eta > 0$, and each choice of $\mathcal{A}_1$ and $\mathcal{A}_2$ with only polynomially many elements,

$$E\Big\{\sup_{x\in\mathcal{A}_1}\sup_{h\in\mathcal{A}_2}|\widehat{F}_h^*(x)-\widetilde{F}(x|h)|\Big\}=O(n^{\eta-(1/2)}).$$

Using this property, and the fact that $K$ is Hölder continuous, it may be shown by a "continuity argument" (see, e.g., [8]) that

$$E\Big\{\sup_{x\in\mathcal{I}}\sup_{h\in[Cn^{-1/5},h_1]}|\widehat{F}_h^*(x)-\widetilde{F}(x|h)|\Big\}=o(1)$$

for any $C > 0$. It follows that if $\tilde{h}$ is a random element of the interval $[Cn^{-1/5}, h_1]$,

(5.31) $$E\Big\{\sup_{x\in\mathcal{I}}|\widehat{F}_{\tilde{h}}^*(x)-\widetilde{F}(x|\tilde{h})|\Big\}=o(1).$$

Write simply $\widehat{F}^*$ for $\widehat{F}^*_{\hat{h}_{\text{crit}}}$, and put $\widehat{H}^* = -\log(1-\widehat{F}^*)$ and $\widetilde{H} = -\log(1-\widetilde{F})$. Taking $h_1 = h_0$ and $\tilde{h} = \hat{h}_{\text{crit}}$, which in view of (5.29) and the assumptions in the theorem satisfies $P(Cn^{-1/5} \leq \hat{h}_{\text{crit}} \leq h_0) \to 1$ for some $C > 0$, we deduce from (5.31) that $|\widehat{F}^* - \widetilde{F}(\cdot|\hat{h}_{\text{crit}})| = o_p(1)$ uniformly on $\mathcal{I}$. From this result and (5.30), we see that

(5.32) $$\sup_{x\in\mathcal{I}}|\widehat{H}^*(x)-\widetilde{H}(x|\hat{h}_{\text{crit}})|=R_1^*,$$

where, here and below, $R_j^*$ denotes a random variable that is defined through Monte Carlo simulation conditional on $\mathcal{X}$ and satisfies $P(|R_j^*| > \eta) \to 0$ for each $\eta > 0$, where the probability is defined in the unconditional sense.

If $T^*$ denotes the bootstrap version of $T$, then

$$T^* = \iint_{x,y\,:\,x+y,x-y\in\mathcal{I}} \max\{0, 2\widehat{H}^*(x)-\widehat{H}^*(x+y)-\widehat{H}^*(x-y)\}\,dx\,dy$$

$$= \iint_{x,y\,:\,x+y,x-y\in\mathcal{I}} \max\{0, 2\widetilde{H}(x|\hat{h}_{\text{crit}})-\widetilde{H}(x+y|\hat{h}_{\text{crit}})$$

$$-\widetilde{H}(x-y|\hat{h}_{\text{crit}})\}\,dx\,dy + R_2^*$$

$$= R_2^*,$$

where the second identity follows from (5.32) and the third from the fact that, by the definition of $\hat{h}_{\text{crit}}$, $\widetilde{H}(\cdot|\hat{h}_{\text{crit}})$ is convex on $\mathcal{I}$. Therefore, $P(T^* > \eta) \to 0$ for each $\eta > 0$. Hence, since $\hat{c}(\alpha)$ is defined by $P\{T^* > \hat{c}(\alpha)|\mathcal{X}\} = \alpha$, then (5.28) must hold.

Centre for Mathematics
  and Its Applications
Australian National University
Canberra, ACT 0200
Australia
e-mail: Peter.Hall@maths.anu.edu.au

Institut de Statistique
Université Catholique de Louvain
Voie du Roman Pays 20
B-1348 Louvain-la-Neuve
Belgium
e-mail: vankeilegom@stat.ucl.ac.be